\documentclass[11pt]{amsart}

\usepackage{amsthm,amssymb,amsfonts,amsmath,enumerate}

\newcommand{\Z}{\mathbb{Z}}
\newcommand{\R}{\mathbb{R}}
\newcommand{\Q}{\mathbb{Q}}
\renewcommand{\phi}{\varphi}
\DeclareMathOperator*{\lcm}{lcm}
\DeclareMathOperator*{\conv}{conv}
\newcommand\Ehr{\operatorname{Ehr}}
\renewcommand{\P}{{\mathcal P}}
\def\th{^{\text{th}}}

\newtheorem{theorem}{Theorem}
\newtheorem{corollary}[theorem]{Corollary}
\newtheorem{proposition}[theorem]{Proposition}
\newtheorem{lemma}[theorem]{Lemma}
\newtheorem{conjecture}[theorem]{Conjecture}

\newtheorem*{definition}{Definition}

\newcommand\comment[1]{}                %  Silent version.
\renewcommand\comment[1]{\emph{[#1]}}           %  Comment revealed.

\title{Maximal Periods of (Ehrhart) Quasi-Polynomials}

\author{Matthias Beck}
\address{Department of Mathematics\\
         San Francisco State University\\
         San Francisco, CA 94132\\
         USA}
\email{beck@math.sfsu.edu}
\urladdr{http://math.sfsu.edu/beck}

\author{Steven V. Sam}
\address{Department of Mathematics\\
  University of California\\
  Berkeley, CA 94720\\
  USA}
\email{ssam@berkeley.edu}

\author{Kevin M. Woods}
\address{Department of Mathematics\\
  Oberlin College\\
  Oberlin, OH 44074\\
  USA}
\email{kevin.woods@oberlin.edu}
\urladdr{http://www.oberlin.edu/math/faculty/woods.html}

\thanks{The authors thank an anonymous referee for helpful suggestions.}


\date{August 13, 2007. To appear in \emph{Journal of Combinatorial Theory Series A}}

\begin{document}

\begin{abstract}
  A \emph{quasi-polynomial} is a function defined of the form $q(k) =
  c_d(k) \, k^d + c_{d-1}(k) \, k^{d-1} + \dots + c_0(k)$, where $c_0,
  c_1, \dots, c_d$ are periodic functions in $k \in \Z$.  Prominent
  examples of quasi-polynomials appear in Ehrhart's theory as
  integer-point counting functions for rational polytopes, and
  McMullen gives upper bounds for the periods of the $c_j(k)$ for
  Ehrhart quasi-polynomials.  For generic polytopes, McMullen's bounds
  seem to be sharp, but sometimes smaller periods exist.  We prove
  that the second leading coefficient of an Ehrhart quasi-polynomial
  always has maximal expected period and present a general theorem
  that yields maximal periods for the coefficients of certain
  quasi-polynomials.  We present a construction for (Ehrhart)
  quasi-polynomials that exhibit maximal period behavior and use it to
  answer a question of Zaslavsky on convolutions of quasi-polynomials.
\end{abstract}

\maketitle

\section{Introduction}

A \emph{quasi-polynomial} is a function defined on $\Z$ of the form
\begin{equation}\label{quasipol}
  q(k) = c_d(k) \, k^d + c_{d-1}(k) \, k^{d-1} + \dots + c_0(k) \, ,
\end{equation}
where $c_0, c_1, \dots, c_d$ are periodic functions in $k$, called the
\emph{coefficient functions} of $q$. Assuming $c_d$ is not identically
zero, we call $d$ the \emph{degree} of~$q$.  Quasi-polynomials play a
prominent role in enumerative combinatorics \cite[Chapter
4]{stanleyec1}.  Arguably their best known appearance is in Ehrhart's
fundamental work on integer-point enumeration in rational polytopes
\cite{ehrhartpolynomial}. For more applications, we refer to the
recent article \cite{lisonekquasipol}.

A \emph{rational polytope} $\P \subset \R^n$ is the convex hull of
finitely many points in $\Q^n$.  The \emph{dimension} of a polytope
$\P$ is the dimension $d$ of the smallest affine space containing
$\P$, in which case we call $\P$ a $d$-polytope. A \emph{face} of $\P$
is a subset of the form $\P \cap H$, where $H$ is a hyperplane such
that $\P$ is entirely contained in one of the two closed half-spaces
of $\R^n$ that $H$ naturally defines. A $(d-1)$-face of a $d$-polytope
is a \emph{facet}, and a $0$-face is a \emph{vertex}. The smallest $k
\in \Z_{ >0 }$ for which the vertices of $k \P$ are in $\Z^n$ is the
\emph{denominator} of $\P$.  Ehrhart's theorem states that the
integer-point counting function $L_\P (k) := \# \left( k \P \cap \Z^n
\right)$ is a quasi-polynomial of degree $d$ in $k \in \Z_{ >0 }$, and
the denominator of $\P$ is a period of each of the coefficient
functions.
For a general introduction to polytopes, we refer to \cite{ziegler};
for an introduction to Ehrhart theory, see \cite{ccd}.

In general, many of the coefficient functions will have smaller
periods. Suppose $q$ is given by \eqref{quasipol}.  The \emph{minimum
  period} of $c_j$ is the smallest $p\in\Z_{>0}$ such that $c_j(k+p) =
c_j(k)$ for all $k \in \Z$ (any multiple of $p$ is, of course, also a
period of $c_j$).  The \emph{minimum period} of $q$ is the least
common multiple of the minimum periods of $c_0, c_1, \dots, c_d$.  In
this paper, we study the minimum periods of the $c_j$. All of our
illustrating examples can be realized as Ehrhart quasi-polynomials.
Ehrhart's theorem tells us that the minimum period of each $c_j$
divides the denominator of $\P$.

The following theorem due to McMullen \cite[Theorem
6]{mcmullenreciprocity} gives a more precise upper bound for these
periods.  For $0\le j\le d$, define the \emph{$j$-index} of $\P$ to be
the minimal positive integer $p_j$ such that the $j$-dimensional faces
of $p_j \P$ all span affine subspaces that contain integer lattice
points.

\begin{theorem}[McMullen] \label{mcmullenthm} Given a rational
  $d$-polytope $\P$, let $p_j$ be the $j$-index of $\P$. If $L_\P (k)
  = c_d(k) \, k^d + c_{d-1}(k) \, k^{d-1} + \dots + c_0(k)$ is the
  Ehrhart quasi-polynomial of $\P$, then the minimum period of $c_j$
  divides~$p_j$.
\end{theorem}

Note that $p_d | p_{ d-1 } | \cdots | p_0$.  Since $p_0$ is the
denominator of $\P$, this is a stronger version of Ehrhart's theorem.
If we further assume that $\P$ is full-dimensional, then $p_d=1$, and
so $c_d(k)$ is a constant function.  In this case, it is well known
that $c_d(k)$ is the Euclidean volume of $\P$
\cite{ccd,ehrhartpolynomial}.

These bounds on the periods seem tight for generic rational polytopes,
that is, $p_j$ is the minimum period of $c_j$, but this statement is
ill-formed (we make no claim what notion of \emph{genericity} should
be used here) and conjectural.  One of the contributions of this paper
is a step in the right direction: for any $p_d | p_{ d-1 } | \cdots |
p_0$, there does indeed exist a polytope such that $c_j$ has minimum
period $p_j$.

\begin{theorem} \label{fullperiodsimplexthm} Given distinct positive
  integers $p_d | p_{ d-1 } | \cdots | p_0$, the simplex
\[
\Delta = \conv \left\{ \left( \tfrac{ 1 }{ p_0 } , 0, \dots, 0 \right)
  , \left( 0, \tfrac{ 1 }{ p_1 } , 0, \dots, 0 \right) , \dots, \left(
    0, \dots, 0, \tfrac{ 1 }{ p_d } \right) \right\}\subset \R^{d+1}
\]
has an Ehrhart quasi-polynomial $L_\Delta (k) = c_d(k) \, k^d +
c_{d-1}(k) \, k^{d-1} + \dots + c_0(k)$, where $c_j$ has minimum
period $p_j$ for $j = 0, 1, \dots, d$ (and $p_j$ is the $j$-index of
$\Delta$).
\end{theorem}
Note that $\Delta$ is actually not a full-dimensional polytope; it is
a $d$-dimensional polytope in $\R^{d+1}$.  This allows us to state the
theorem in slightly greater generality (we don't have to constrain
$p_d=1$, which is necessary for a full-dimensional polytope).

Theorem \ref{fullperiodsimplexthm} complements recent literature
\cite{deloeramcallister, mcallisterwoods} that contains several
special classes of polytopes that defy the expectation that $c_j$ has
minimum period $p_j$.  De Loera--McAllister \cite{deloeramcallister}
constructed a family of polytopes stemming from representation theory
that exhibit \emph{period collapse}, i.e., the Ehrhart
quasi-polynomials of these polytopes (which have arbitrarily large
denominator) have minimum period 1---they are polynomials.
McAllister--Woods \cite{mcallisterwoods} gave a class of polytopes
whose Ehrhart quasi-polynomials have arbitrary period collapse (though
not for the periods of the individual coefficient functions), as well
as an example of non-monotonic minimum periods of the coefficient
functions.

First, we will prove (in Section \ref{secondcoeffsection}) that no
period collapse is possible in the second leading coefficient
$c_{d-1}(k)$:

\begin{theorem} \label{secondcoeffthm} Given a rational $d$-polytope
  $\P$, let $p_{d-1}$ be the $(d-1)$-index of $\P$. Let $L_\P (k) =
  c_d(k) \, k^d + c_{d-1}(k) \, k^{d-1} + \dots + c_0(k)$. Then $c_{
    d-1 }$ has minimum period $p_{d-1}$.
\end{theorem}

In Section \ref{generalquasithmsection}, we give some general results
on quasi-polynomials with maximal period behavior. Namely, we will
prove:

\begin{theorem}\label{quasimonomialthm}
  Suppose $c(k)$ is a periodic function with minimum period $n$, and
  $m$ is some nonnegative integer. Then the rational generating
  function $\sum_{k\geq 0} c(k)k^mx^k$ has as poles only $n\th$ roots
  of unity, and each of these poles has order $m+1$.
\end{theorem}

A direct consequence of this statement is the following:

\begin{corollary}\label{generalquasithm}
  Suppose $r(x)$ is a proper rational function all of whose poles are
  primitive $n\th$ roots of unity. Then $r$ is the generating function
  of a quasi-polynomial
\[
r(x) = \sum_{ k \ge 0 } \left( c_d(k) \, k^d + c_{d-1}(k) \, k^{d-1} +
  \dots + c_0(k) \right) x^k ,
\]
where each $c_j$ is either identically zero or has minimum period $n$.
\end{corollary}

As an application to Theorem~\ref{fullperiodsimplexthm} (proved in
Section \ref{fullperiodsimplexsection}), we turn to a question that
stems from a recent theorem of Zaslavsky
\cite{zaslavskyquasipolynomial}.  Suppose $A(k) = a_d(k) \, k^d +
a_{d-1}(k) \, k^{d-1} + \dots + a_0(k)$ and $B(k) = b_e(k) \, k^e +
b_{e-1}(k) \, k^{e-1} + \dots + b_0(k)$ are quasi-polynomials, where
the minimum period of $a_j$ is $\alpha_j$ and the minimum period of
$b_j$ is $\beta_j$. Then the \emph{convolution}
\[
  C(k) := \sum_{ m=0 }^k A(k-m) \, B(m)
\]
is another quasi-polynomial. If we write $C(k) = c_{d+e+1}(k) \,
k^{d+e+1} + c_{d+e}(k) \, k^{d+e} + \dots + c_0(k)$, and let $c_j$
have minimum period $\gamma_j$, Zaslavsky proved the following result.

\begin{theorem}[Zaslavsky] \label{zaslavthm} Define
  $g_j=\lcm\{\gcd(\alpha_i, \beta_{j-i}) : \, 0\le i \le d,\ 0\le j-i
  \le e\}$ for $j \ge 0$, and let $g_{-1}=1$. Then
\begin{equation}\label{zasdivide}
  \gamma_{j+1} \left| \, \lcm
    \left\{\alpha_{j+1},\dots,\alpha_d,\beta_{j+1}, \dots, \beta_e,g_j
    \right\} \right. .
\end{equation}
\end{theorem}

We will reprove this result in Section \ref{zaslavskysection} using
the generating-function tools we develop.  A natural problem, raised
by Zaslavsky, is to construct two quasi-polynomials whose convolution
satisfies \eqref{zasdivide} with equality. The answer is given by
another application of Theorem~\ref{fullperiodsimplexthm} (Section
\ref{zaslavskysection}).

\begin{theorem} \label{zaslavsharpthm} Given $d \ge e$ and distinct
  positive integers $\alpha_d | \alpha_{ d-1 } | \cdots | \alpha_e |
  \beta_e | \alpha_{ e-1 } | \beta_{ e-1 } | \cdots | \alpha_0 |
  \beta_0$, let
\[
\Delta_1 = \conv \left\{ \left( \tfrac{ 1 }{ \alpha_0 } , 0, \dots, 0
  \right) , \left( 0, \tfrac{ 1 }{ \alpha_1 } , 0, \dots, 0 \right) ,
  \dots, \left( 0, \dots, 0, \tfrac{ 1 }{ \alpha_d } \right) \right\}
\]
and
\[
\Delta_2 = \conv \left\{ \left( \tfrac{ 1 }{ \beta_0 } , 0, \dots, 0
  \right) , \left( 0, \tfrac{ 1 }{ \beta_1 } , 0, \dots, 0 \right) ,
  \dots, \left( 0, \dots, 0, \tfrac{ 1 }{ \beta_e } \right) \right\} .
\]
Then the convolution of $L_{ \Delta_1 }$ and $L_{ \Delta_2 }$
satisfies \eqref{zasdivide} with equality.
\end{theorem}

%%%%%%%%%%%%%%%%%%%%%%%%%%%%%%%%%%%%%%%%%%%%%%%%%%%%

\section{The Second Leading Coefficient of an Ehrhart
  Quasi-Polynomial}\label{secondcoeffsection}

In this section we prove Theorem~\ref{secondcoeffthm}, namely the
minimum period of the second leading coefficient of the Ehrhart
quasi-polynomial of a rational $d$-polytope $\P$ equals the
$(d-1)$-index of $\P$.
% (the minimal positive integer $p_{d-1}$ such that the facets of
% $p_{d-1} \P$ all span affine subspaces that contain lattice points).
Most of the work towards Theorem~\ref{secondcoeffthm}
is contained in the proof of the following result.

\begin{proposition} \label{secondcoeffprop} If $\P$ is a rational
  $d$-polytope with Ehrhart quasi-polynomial $L_\P (k) = c_d(k) \, k^d
  + c_{d-1}(k) \, k^{d-1} + \dots + c_0(k)$, then $c_{ d-1 }$ is
  constant if and only if the $(d-1)$-index of $\P$ is $1$.
\end{proposition}

\begin{proof}
  If the $(d-1)$-index of $\P$ is $1$, then $c_{ d-1 }$ is constant by
  McMullen's Theorem~\ref{mcmullenthm}.

  For the converse implication, we use the \emph{Ehrhart-Macdonald
    Reciprocity Theorem} \cite{ccd,macdonald}. It says that for a
  rational $d$-polytope $\P$, the evaluation of $L_\P$ at negative
  integers yields the lattice-point enumerator of the interior
  $\P^\circ$, namely,
\[
  L_\P (-k) = (-1)^d L_{ \P^\circ } (k) \, .
\]
This identity implies that the lattice-point enumerator for the
boundary of $\P$ is the quasi-polynomial $L_{ \partial \P } (k) = L_\P
(k) - (-1)^d L_\P (-k)$. Since $L_{\partial \P}(k)$ counts integer
points in a $(d-1)$-dimensional object, it is a degree $d-1$
quasi-polynomial, and we see that its leading coefficient is $c_{ d-1
}(k)+c_{d-1}(-k)$.

%\comment{Referee was confused by our terminology ``leading
%  coefficient'' in what follows.}

%Now suppose the $(d-1)$-index of $\P$ is $m>1$. Then there is a facet
%of $\P$ whose affine span has no lattice points when dilated by $1,
%m+1, 2m+1, \dots$ On the other hand, the affine span of \emph{any}
%facet contains lattice points when dilated by multiples of $m$.  The
%leading coefficient of the quasi-polynomial $L_{ \partial \P } (mk)$
%is a constant measuring the volume of the facets relative to the
%sublattices in their affine spans (for the same reason that the
%leading coefficient of $L_{\P}(k)$ measures the volume of $\P$:
%asymptotically, the volume is a good approximation for the number of
%integer points).  The same can be said for the leading coefficient of
%the polynomial $L_{ \partial \P } (1+mk)$; however, now some of these
%sublattices are empty, and thus the leading coefficient of
%$L_{ \partial \P } (1+mk)$ is smaller than the leading coefficient of
%$L_{ \partial \P } (mk)$, i.e., $c_{ d-1 } (1+mk) + c_{d-1}(-1-mk)<\,
%c_{ d-1 } (mk) + c_{d-1}(-mk)$. Hence $c_{ d-1 }$ is not constant.

Suppose that the $(d-1)$-index of $\P$ is $m>1$, and that $c_{d-1}$ is
a constant. Then the leading coefficient of $L_{ \partial \P }(k)$ is
constant, and the affine span of every facet of $\P$ contains lattice
points when dilated by any multiple of $m$. However, there are facets
of $\P$ whose affine spans contain no lattice points when dilated by
$jm+1$ for $j \ge 0$. Let $F_1, \dots, F_n$ be these facets, and
consider the polytopal complex $\P' = \bigcup F_i$. In fact, the
lattice points of $k\P' := \bigcup kF_i$ are counted by a
quasi-polynomial $L_{\P'}(k)$. We can obtain $L_{\P'}(k)$ by first
starting with $L_{ \partial \P }(k)$. Then for each facet of $\P$ not
among $F_1, \dots, F_n$, subtract its Ehrhart quasi-polynomial from
$L_{ \partial \P }(k)$. Some of the lower dimensional faces of $\P'$
might now be uncounted by the resulting enumerator, so we play an
inclusion-exclusion game with their Ehrhart quasi-polynomials to get
$L_{\P'}(k)$ as a sum of Ehrhart quasi-polynomials of the faces of
$\P$. We are concerned only with the leading coefficient function of
$L_{\P'}(k)$, which is unaffected by this inclusion-exclusion. The
Ehrhart quasi-polynomial for each facet not among $F_1, \dots, F_n$
has constant leading term by McMullen's Theorem, so the leading term
of $L_{\P'}(k)$ is some constant $c$. This means that for large values
of $k$, the number of lattice points in $k\P'$ is asymptotically
$c \, k^{d-1}$. However, by construction of $\P'$, we have $L_{\P'}(jm +
1) = 0$ for all $j \ge 0$, which gives a contradiction. Thus, if the
$(d-1)$-index of $\P$ is greater than 1, then $c_{d-1}$ is not a
constant.
\end{proof}

\begin{proof}[Proof of Theorem \ref{secondcoeffthm}]
  Let $p$ be the minimal period of $c_{ d-1 }$ and $q$ be the
  $(d-1)$-index of $\P$.  By McMullen's Theorem \ref{mcmullenthm},
  $p|q$.  On the other hand, the second-leading coefficient of $L_{ p
    \P }$ is constant, and by Proposition~\ref{secondcoeffprop}, the
  $(d-1)$-index of $p \P$ is 1, which implies $q|p$.
\end{proof}

%%%%%%%%%%%%%%%%%%%%%%%%%%%%%%%%%%%%%%%%%%%%%%%%%%%%

\section{Some General Results on Quasi-Polynomial
  Periods}\label{generalquasithmsection}

A key ingredient to proving Theorem~\ref{quasimonomialthm} is a basic
result (see, e.g., \cite[Chapter 3]{ccd} or \cite[Chapter
4]{stanleyec1}) about a quasi-polynomial $q(k)$ and its generating
function $r(x)=\sum_{k\ge 0}q(k)x^k$, which is easily seen to be a
rational function.

\begin{lemma} \label{basiclemma} Suppose $q$ is a quasi-polynomial
  with generating function $r(x) = \sum_{ k \ge 0 } q(k) \,
  x^k$ (which evaluates to a proper rational function).
  Then $n$ is a period of $q$ and $q$ has degree $d$ if and only
  if all poles of $r$ are $n\th$ roots of unity of order $\le d+1$ and
  there is a pole of order $d+1$.
\end{lemma}

The above result will be useful again in the proof of
Theorem~\ref{fullperiodsimplexthm}. Recall that the statement of
Theorem~\ref{quasimonomialthm} is that given a periodic function $c(k)$
with minimum period $n$ and a nonnegative integer $m$, the only poles
of the rational generating function $\sum_{k\geq 0} c(k)k^mx^k$
are $n\th$ roots of unity, and each pole has order $m+1$.

\begin{proof}[Proof of Theorem~\ref{quasimonomialthm}]
  We use induction on $m$. The case $m=0$ follows directly from Lemma
  \ref{basiclemma}, as
  \[ \sum_{ k \ge 0 } c(k) k^0
  x^k=\frac{c(0)+c(1)x+\cdots+c(n-1)x^{n-1}}{1-x^n} \, .\]
  The induction step is a consequence of the identity
  \[
  \sum_{ k \ge 0 } c(k) k^m x^k = x \, \frac{ d }{ dx } \sum_{ k \ge 0
  } c(k) k^{m-1} x^k
  \]
  and the fact that a pole of order $m-1$ turns into a pole of order
  $m$ under differentiation.
\end{proof}

Corollary~\ref{generalquasithm} now follows like a breeze. Recall its
statement: If $r(x)$ is a proper rational function all of whose poles are
primitive $n\th$ roots of unity, then $r$ is the generating function
of a quasi-polynomial
\[
r(x) = \sum_{ k \ge 0 } \left( c_d(k) \, k^d + c_{d-1}(k) \, k^{d-1} +
  \dots + c_0(k) \right) x^k ,
\]
where each $c_j \not\equiv 0$ has minimum period $n$.

\begin{proof}[Proof of Corollary \ref{generalquasithm}]
  % Lemma \ref{basiclemma} implies that $p|n$.
  Consider the rational generating functions
\[
r_j (x) := \sum_{ k \ge 0 } c_j(k) k^j x^k , \qquad \text{ so that }
\qquad r(x) = r_d(x) + r_{ d-1 }(x) + \dots + r_0(x) \, .
\]
We claim that the poles of each (not identically zero) $r_j(x)$ are
all primitive
$n\th$ roots of unity. Indeed, suppose not, and consider the largest
$j$ such that $r_j(x)$ has a pole $\omega$ which is not a primitive
$n\th$ root of unity. Theorem~\ref{quasimonomialthm} says that
$\omega$ is a pole of $r_j(x)$ of order $j+1$.  Since $\omega$ is not
a pole of $r_d(x),r_{d-1}(x),\ldots,r_{j+1}(x)$ (we chose $j$ as large
as possible), $\omega$ is a pole of
\[r_d(x)+r_{d-1}(x)+\cdots+r_{j+1}(x)+r_j(x)\] of order $j+1$.  On the
other hand, Theorem~\ref{quasimonomialthm} also implies that
$r_{j-1}(x),r_{j-2}(x),\ldots,r_0(x)$ have no poles of order greater
than $j$.  Summing over all the $r_i$, $\omega$ must be a pole of
$r(x)$ of order $j+1$, contradicting the fact that $r(x)$ has only
poles that are primitive
$n\th$ roots of unity.

Therefore the poles of each (not identically zero) $r_j(x)$ are all
primitive roots of unity.  Lemma~\ref{basiclemma} implies that $n$ is
a period of each nonzero $c_j$.
\end{proof}

%%%%%%%%%%%%%%%%%%%%%%%%%%%%%%%%%%%%%%%%%%%%%%%%%%%%

\section{Ehrhart Quasi-Polynomials with Maximal
  Periods}\label{fullperiodsimplexsection}

Recall that Theorem~\ref{fullperiodsimplexthm} says that for given
distinct positive integers $p_d | p_{ d-1 } | \cdots | p_0$, the
simplex
\[
\Delta = \conv \left\{ \left( \tfrac{ 1 }{ p_0 } , 0, \dots, 0 \right)
  , \left( 0, \tfrac{ 1 }{ p_1 } , 0, \dots, 0 \right) , \dots, \left(
    0, \dots, 0, \tfrac{ 1 }{ p_d } \right) \right\}\subset \R^{d+1}
\]
has the Ehrhart quasi-polynomial $L_\Delta (k) = c_d(k) \, k^d +
c_{d-1}(k) \, k^{d-1} + \dots + c_0(k)$, where $c_j$ has minimum
period $p_j$ for $j = 0, 1, \dots, d$. Note that $p_j$ is the
$j$-index of $\Delta$.

\begin{proof}[Proof of Theorem \ref{fullperiodsimplexthm}]
The Ehrhart series of
\[
\Delta = \left\{ \left( x_0, x_1, \dots, x_d \right) \in \R_{ \ge 0
  }^{ d+1 } : \, p_0 x_0 + p_1 x_1 + \dots + p_d x_d = 1 \right\}
\]
is, by construction,
\[
\Ehr_\Delta(x) := \sum_{ k \ge 0 } L_\Delta(k) \, x^k = \frac{ 1 }{
  \left( 1 - x^{ p_0 } \right) \left( 1 - x^{ p_1 } \right) \cdots
  \left( 1 - x^{ p_d } \right) } \, .
\]
Given $j$, let $\omega$ be a primitive $p_j\th$ root of unity. Then
$\omega$ is a pole of $\Ehr_\Delta(x)$ of order $j+1$. We expand
$\Ehr_\Delta (x)$ to yield the Ehrhart quasi-polynomial:
\[
\Ehr_\Delta(x) = \sum_{ k \ge 0 } L_\Delta(k) \, x^k = \sum_{ k \ge 0
} \left( c_d(k) \, k^d + c_{d-1}(k) \, k^{d-1} + \dots + c_0(k)
\right) x^k .
\]

Let $n$ be the minimum period of $c_j(k)$.  By McMullen's Theorem
\ref{mcmullenthm}, $n|p_j$.  Therefore, we need to show that $p_j|n$.
As before, let $r_j(x)=\sum_{k\ge 0}c_j(k)k^jx^k$, so that
$\Ehr_{\Delta}(x)=r_d(x)+r_{d-1}(x)+\cdots+r_0(x).$ Since $\omega$ is
a pole of $\Ehr_{\Delta}(x)$, it must be a pole of (at least) one of
$r_d,\ldots,r_0$.  Let $J$ be the largest index such that $\omega$ is
a pole of $r_J(x)$. By Theorem~\ref{quasimonomialthm}, $\omega$ is a
pole of $r_J(x)$ of order $J+1$.  Since $\omega$ is not a pole of
$r_d(x),r_{d-1}(x),\ldots,r_{J+1}(x)$, $\omega$ is a pole of
\[r_d(x)+r_{d-1}(x)+\cdots+r_{J+1}(x)+r_J(x)\] of order $J+1$.  On the
other hand, Theorem~\ref{quasimonomialthm} also implies that
$r_{J-1}(x),r_{J-2}(x),\ldots,r_0(x)$ have no poles of order greater
than $J$.  Summing over all the $r_i$, $\omega$ must be a pole of
$\Ehr_\Delta(x)$ of order $J+1$.  Since we saw that $\omega$ is a pole
of $\Ehr_\Delta(x)$ of order $j+1$, we have that $J=j$, that is,
$\omega$ is a pole of $r_j(x)$.  Since $\omega$ is a primitive
$p_j\th$ root of unity, Theorem~\ref{quasimonomialthm} says that $p_j$
must divide the minimum period $n$, and so $n=p_j$, as desired.
\end{proof}

%%%%%%%%%%%%%%%%%%%%%%%%%%%%%%%%%%%%%%%%%%%%%%%%%%%%

\section{Quasi-Polynomial Convolution with Maximal
  Periods}\label{zaslavskysection}

We start our last section with a generating-function proof of
Zaslavsky's Theorem \ref{zaslavthm}. It uses the following
generalization of Lemma \ref{basiclemma}:

\begin{lemma} \label{basiclemmarefined} Suppose $q(k) = c_d(k) \, k^d
  + c_{ d-1 }(k) \, k^{ d-1 } + \dots + c_0(k)$ is a quasi-polynomial
  with rational generating function $r(x) = \sum_{ k \ge 0 } q(k) \,
  x^k$.
\begin{enumerate}[{\rm (a)}]
\item \label{refineda} If $n$ is a period of $c_j$, then there is an
  $n\th$ root of unity that is a pole of $r$ of order at least $j+1$.
\item \label{refinedb} If all poles of $r$ of order $\ge j+1$ are
  $n\th$ roots of unity, then $n$ is a period of $c_j$.
\end{enumerate}
\end{lemma}

\begin{proof}
Part \eqref{refineda} follows from Theorem~\ref{quasimonomialthm}.

For part \eqref{refinedb}, expand $r$ (crudely) into partial fractions
as $r(x) = s(x) + t(x)$, such that $s$ has as poles the poles of $r$
of order $\ge j+1$ and $t$ has as poles those of order $\le j$. Now
apply Lemma \ref{basiclemma} to $s$ and note that $t$ does not
contribute to $c_j$.
\end{proof}

\begin{proof}[Proof of Theorem \ref{zaslavthm}]
  Let $f_A(x) = \sum_{ k \ge 0 } A(k) \, x^k$ and define $f_B$ and
  $f_C$ analogously. To determine $\gamma_{ j+1 }$, the period of $c_{
    j+1 }$, Lemma \ref{basiclemmarefined}\eqref{refinedb} tells us
  that we need to consider the poles of $f_C(x) = f_A(x) f_B(x)$ of
  order $\ge j+2$. These poles come in three types:
\begin{enumerate}[{\rm (1)}]
\item poles of $f_A$ of order $\ge j+2$;
\item poles of $f_B$ of order $\ge j+2$;
\item common poles of $f_A$ and $f_B$ whose orders add up to at least
  $j+2$.
\end{enumerate}
Lemma~\ref{basiclemmarefined}\eqref{refineda} gives the statement of
Theorem~\ref{zaslavthm} instantly; the periods
$\alpha_{j+1},\dots,\alpha_d$ give rise to poles of type (1),
$\beta_{j+1},\dots,\beta_e$ give rise to poles of type (2), and
$g_j=\lcm\{\gcd(\alpha_i, \beta_{j-i}) : \, 0\le i \le d,\ 0\le j-i
\le e\}$ stems from poles of type (3).
\end{proof}

\begin{proof}[Proof of Theorem \ref{zaslavsharpthm}]
  The convolution of $L_{ \Delta_1 }$ and $L_{ \Delta_2 }$ equals
  $L_\Delta$, where $\Delta$ is the $(d+e+1)$-simplex
\[
\Delta = \conv \left\{ \left( \tfrac{ 1 }{ \alpha_0 } , 0, \dots, 0
  \right) , \dots, \left( 0, \dots, 0, \tfrac{ 1 }{ \alpha_d }, 0,
    \dots, 0 \right) , \left( 0, \dots, 0, \tfrac{ 1 }{ \beta_0 } , 0,
    \dots, 0 \right) , \dots, \left( 0, \dots, 0, \tfrac{ 1 }{ \beta_e
    } \right) \right\} ,
\]
which follows directly from the fact that the generating function of
the convolution of two quasi-polynomials is the product of their
generating functions.  Let
\[
  L_\Delta(k) = c_{ d+e+1 }(k) \, k^{ d+e+1 } + c_{ d+e }(k) \, k^{ d+e } + \dots + c_0(k)
\]
and suppose $c_j(k)$
has minimum period $\gamma_j$.  By construction and Theorem
\ref{fullperiodsimplexthm}, we have
\[
\gamma_{ 2j } = \beta_j \qquad \text{ and } \qquad \gamma_{ 2j+1 } =
\alpha_j \qquad \text{ for } 0 \le j \le e \, ,
\]
and $\gamma_{ e+j+1 } = \alpha_j$ for $j>e$.  We will show that these
values agree with the upper bounds given by Zaslavsky's Theorem
\ref{zaslavthm}.  We distinguish three cases.

Case 1: $j \le 2e$ and $j+1 = 2m$ for some integer $m$.
We need to show that
\begin{equation} \label{firstlcm} \gamma_{ j+1 } = \lcm \left\{
    \alpha_{ j+1 } , \alpha_{ j+2 } , \dots, \alpha_d, \beta_{ j+1 } ,
    \beta_{ j+2 } , \dots, \beta_e, g_j \right\} = \beta_m \, .
\end{equation}
Consider
\[
g_j = \lcm \left\{ \gcd \left( \alpha_i , \beta_{ j-i } \right) : \, 0
  \le i \le d, \, 0 \le j-i \le e \right\} .
\]
If $2i \ge j$, i.e., $i \ge m$, then $\gcd \left( \alpha_i , \beta_{
    j-i } \right) = \beta_{ j-i }$. Thus
\[
g_j = \lcm \left\{ \alpha_j, \alpha_{ j-1 }, \dots, \alpha_{ m+1 },
  \beta_m, \beta_{ m+1 }, \dots, \beta_j \right\} = \beta_m \, ,
\]
which proves \eqref{firstlcm}, since $j+1 > m$.

Case 2: $j \le 2e$ and $j = 2m$ for some integer $m$.
We need to show that
\begin{equation} \label{secondlcm} \gamma_{ j+1 } = \lcm \left\{
    \alpha_{ j+1 } , \alpha_{ j+2 } , \dots, \alpha_d, \beta_{ j+1 } ,
    \beta_{ j+2 } , \dots, \beta_e, g_j \right\} = \alpha_m \, .
\end{equation}
Now
\[
g_j = \lcm \left\{ \alpha_j, \alpha_{ j-1 }, \dots, \alpha_m, \beta_{
    m+1 }, \beta_{ m+2 }, \dots, \beta_j \right\} = \alpha_m \, ,
\]
which proves \eqref{secondlcm}, since $j+1 > m$.

Case 3: $j > 2e$.
We would like to show that
\begin{equation} \label{thirdlcm} \gamma_{ j+1 } = \lcm \left\{
    \alpha_{ j+1 } , \alpha_{ j+2 } , \dots, \alpha_d, \beta_{ j+1 } ,
    \beta_{ j+2 } , \dots, \beta_e, g_j \right\} = \alpha_{ j-e } \, .
\end{equation}
Here
\[
g_j = \lcm \left\{ \gcd \left( \alpha_i , \beta_{ j-i } \right) : \,
  j-e \le i \le j \right\} .
\]
However, for $j-e \le i \le j$, we have $\gcd \left( a_i, \beta_{ j-i
  } \right) = \alpha_i$, whence $g_j = \alpha_{ j-e }$, which proves
\eqref{thirdlcm}.
\end{proof}

%%%%%%%%%%%%%%%%%%%%%%%%%%%%%%%%%%%%%%%%%%%%%%%%%%%%

\section{Open Problems}

For an Ehrhart quasi-polynomial, period collapse cannot happen in
relation to the $j$-index for the first two coefficients. On the other
side, McAllister--Woods \cite{mcallisterwoods} showed that period
collapse can happen for any other coefficient, however, it is still a
mystery to us to what extent.  Tyrrell McAllister \cite{tyrell}
constructed polygons whose Ehrhart periods are $(1, s, t)$ (the
minimum periods of $c_2(k)$, $c_1(k)$, and $c_0(k)$, respectively).

In constructing the simplex with maximal period behavior, we required
that the integers $p_0,\dots,p_d$ be distinct, but perhaps this
restriction is not necessary. Does the statement still hold true if we
weaken the conditions, or do there exist counterexamples?

In the example of periods of quasi-polynomial convolution, Theorem
\ref{zaslavsharpthm}, our methods require that we assume that
$\alpha_d | \alpha_{ d-1 } | \cdots | \alpha_e | \beta_e | \alpha_{
  e-1 } | \beta_{ e-1 } | \cdots | \alpha_0 | \beta_0$, rather than
the more natural $\alpha_d | \alpha_{ d-1 } | \cdots | \alpha_0$ and
$\beta_e | \beta_{ e-1 } | \cdots | \beta_0$.  We conjecture that the
theorem is still true in this case.

More generally, this would follow from a conjecture about a special
class of generating functions:

\begin{conjecture}
  Let $a_1,a_2,\ldots,a_n$ be given positive integers.  Let $q(k) =
  c_d(k) \, k^d + \dots + c_0(k)$ be the
  quasi-polynomial whose generating function $r(x)=\sum_{k\ge
    0}q(k)x^k$ is given by
\[\frac{1}{(1-x^{a_1})(1-x^{a_2})\cdots (1-x^{a_n})} \, . \]
For a positive integer $m$, define $b_m=\#\{i:\ m\mid a_i\}$.  For
$0\le j\le d$, let $p_j=\lcm\{m:\ b_m>j\}$.  Then the minimum period
of $c_j(k)$ is $p_j$.
\end{conjecture}

There are several multi-parameter versions of Ehrhart
quasi-polynomials to which a generalization of McMullen's
Theorem~\ref{mcmullenthm} applies (see \cite[Theorem
7]{mcmullenreciprocity} and \cite{sturmfelsvectorpartitionfunction}).  Beyond
McMullen's theorem, not much is known about periods and minimum
periods (which are now lattices in some $\Z^m$) of these multivariate
quasi-polynomials and coefficient functions.

%%%%%%%%%%%%%%%%%%%%%%%%%%%%%%%%%%%%%%%%%%%%%%%%%%%%

\bibliographystyle{amsplain}
%\bibliography{bib}

\def\cprime{$'$} \def\cprime{$'$}
\providecommand{\bysame}{\leavevmode\hbox to3em{\hrulefill}\thinspace}
\providecommand{\MR}{\relax\ifhmode\unskip\space\fi MR }
\providecommand{\MRhref}[2]{%
  \href{http://www.ams.org/mathscinet-getitem?mr=#1}{#2}
}
\providecommand{\href}[2]{#2}

\end{document}